\documentclass[12pt,fleqn, a4]{article}

\usepackage{amsfonts}
\usepackage{stmaryrd}
\usepackage{amssymb}
\usepackage{latexsym,amscd}
\usepackage{color}
\usepackage{dsfont}
\usepackage{amsbsy}
\usepackage{latexsym}
\usepackage{amssymb}
\usepackage{latexsym}

\usepackage{mathrsfs}
\usepackage{bbm}
\usepackage{textcomp}
\usepackage{mathabx}

\usepackage{graphics}
\usepackage{color}
\newcommand{\colval}{0.3}
\definecolor{colone}{gray}{\colval}

\newcommand{\dcb}{\begin{array}{lll}}
\newcommand{\dce}{\end{array}}
\newcommand{\ebe}{\begin{enumerate}\setlength{\baselineskip}{13pt}\setlength{\parskip}{11pt}}
\newcommand{\dbe}{\end{enumerate}}

\newcommand{\ibegin}{\begin{itemize}\setlength{\baselineskip}{19pt}\setlength{\parskip}{7pt}}
\newcommand{\iend}{\end{itemize}}
\newcommand{\ok}{\rule{4pt}{6pt}}
\newcommand{\desb}{\begin{description}}
\newcommand{\dese}{\end{description}}

\newtheorem{Theorem}{Theorem}[section]
\newtheorem {Cor}[Theorem]{Corollary}
\newtheorem {definition}[Theorem]{Definition}
\newtheorem {pro}[Theorem]{Proposition}
\newtheorem {Lemma}[Theorem]{Lemma}
\newtheorem {rem}[Theorem]{Remark}
\newtheorem {assumption}[Theorem]{Assumption}
\newtheorem {condition}[Theorem]{Condition}
\newtheorem {result}[Theorem]{Result}

\newcommand {\bd}{\begin{definition}}
\newcommand {\ed}{\end{definition}}
\newcommand {\bpro}{\begin{pro}}
\newcommand {\epro}{\end{pro}}
\newcommand {\bl}{\begin{Lemma}}
\newcommand {\el}{\end{Lemma}}
\newcommand {\bcor}{\begin{Cor}}
\newcommand {\ecor}{\end{Cor}}
\newcommand {\brem }{\begin{rem} \rm }
\newcommand {\erem }{\end{rem}}
\newcommand{\bethe}{\begin{Theorem}}
\newcommand{\ethe}{\end{Theorem}}
\newcommand {\bassumption}{\begin{assumption}}
\newcommand {\eassumption}{\end{assumption}}
\newcommand {\bcond}{\begin{condition}}
\newcommand {\econd}{\end{condition}}
\newcommand {\brt}{\begin{result}}
\newcommand {\ert}{\end{result}}

\def \ind{1\!\!1}

\def\lri{[\![}
\def\rri{]\!]}

\begin{document}

\textbf{\Large From Doob's maximal identity to Azema supermartingale}\\

\textbf{Shiqi Song}\\
{\footnotesize Laboratoire de Math\'ematiques et Mod\'elisation d'\'Evry, 91037 \'Evry Cedex, France}

\

\textbf{Abstract}
{\small
Let $\rho$ be a random time. With respect to a filtration $\mathbb{F}$, we study the conditions under which the Azema supermartingale $Z^\rho:=\mathbb{P}[t<\rho\ |\mathcal{F}_{t}], t\geq 0,$ can be written in the form $Z^\rho=\frac{U}{U^*}$, where $U$ is a non negative local martingale starting from 1 vanishing at infinity and $U^*_{t}=\sup_{s\leq t}U_{s}, t\geq 0$. 
}

\



\

\section{Introduction}

Concerning a non negative local martingale $U$, we have the famous Doob maximal inequality about the running maximum process $U^*_{t}=\sup_{s\leq t}U_{s}, t\geq 0$ (cf. \cite[Theorem 2.49]{HWY}, \cite[Theorem 3.8]{KS}, \cite[Theorem 1.7]{RY}). Doob's maximal inequality involves in fact several different inequalities. The most fundamental one is an estimation of the probabilities $\mathbb{P}[U^*_{\infty}>\lambda], \lambda>0$. When $U$ is continuous and vanishes at infinity, applying the Doob optional sampling theorem, we obtain an identity $\mathbb{P}[\sup_{t<s}U_{s}\geq \lambda\ |\mathcal{F}_{t}]=\frac{U_{t}}{\lambda}$ on the set $\{\lambda\geq U_{t}\}$. See \cite[Chapiter 1, Problem 3.28]{KS}. Actually, as indicated in \cite{NY} and \cite{K2015}, this last identity holds for all non negative local martingale $U$ vanishes at infinity, whenever (if and only if) its running maximum process $U^*$ is continuous. 

In \cite[Lemma 2.1]{NY} the above identity of $\mathbb{P}[\sup_{t<s}U_{s}\geq \lambda\ |\mathcal{F}_{t}]$ is linked to the Azema supermartingale $\mathbb{P}[t<\rho\ |\mathcal{F}_{t}]$ of the last maximum time $\rho=\sup\{U=U^*\}$. Actually, $\mathbb{P}[t<\rho\ |\mathcal{F}_{t}]=\mathbb{P}[\sup_{t<s}U_{s}\geq U^*_{t}\ |\mathcal{F}_{t}]=\frac{U_{t}}{U^*_{t}}$. That identity, named Doob maximal identity since \cite{NY}, has been exploited fruitfully in \cite{AP, K2015, NY} to study the last maximum time $\rho$ and the associated filtration enlargement.

Another way that the Doob maximal identity has profited the study of random times  (i.e., random variables taking values in $[0,\infty]$) and filtrations, is that it writes an Azema supermartingale into the multiplication of a local martingale and a non increasing process (a multiplicative decomposition). The multiplicative decomposition of Azema supermartingale, as well as its Doob-Meyer additive decomposition, are among the most fundamental properties of a random time. Various multiplicative decompositions of Azema supermartingale exist, from predictable one (cf. \cite{jacod, PR, song-1405}) to optional one (cf. \cite{K2010nu}). The Doob maximal identity obtained in \cite{AP, K2015, NY} coincides with the predictable multiplicative decomposition. Multiplicative decompositions are all composed of two factors. The Doob maximal identity is distinguished from the others in its direct writing of the second factor as a functional of the first one. This is a very useful property for computations. See \cite{FJS, NY} for some applications.

The Doob maximal identity will not hold without the continuity of the running maximum process. We are interested in whether the identity $\mathbb{P}[t<\rho\ |\mathcal{F}_{t}]=\frac{U_{t}}{U^*_{t}}$ can, nevertheless, remain valid, for a pair of a random time $\rho$ and a local martingale $U$ (which are not necessarily linked via the last maximum time). Precisely, we introduce the following condition for a general random time $\rho$.

\bcond(\textbf{Martingale maximal representation})\label{MMR}
The Azema supermartingale of $\rho$ satisfies $\mathbb{P}[t<\rho\ |\mathcal{F}_{t}]=\frac{U_{t}}{U^*_{t}}, t\geq 0$, for some non negative local martingale $U$ starting from 1 vanishing at infinity. 
\econd

We now want to know what are the random times which satisfy this martingale maximal representation.

\

\section{ABC of random times}\label{ABC}

We work on a stochastic basis $(\Omega,\mathcal{A}, \mathbb{P}, \mathbb{F})$, where $(\Omega,\mathcal{A})$ is a measurable space, $\mathbb{P}$ is a probability measure on the $\sigma$-algebra $\mathcal{A}$, and $\mathbb{F}=(\mathcal{F}_{t})_{t\in\mathbb{R}_{+}}$ is a filtration of sub-$\sigma$-algebras of $\mathcal{A}$. We suppose the usual conditions for $\mathbb{F}$ and we set $\mathcal{F}_{\infty}=\sigma(\cup_{t\geq 0}\mathcal{F}_{t})$. We call a random variable in $[0,\infty]$ a random time. We apply the theory of semimartingale calculus as presented in \cite{DM2,HWY,jacod,protter}.

\subsection{The four projections}\label{basic-elements}

Since the pioneering works such as \cite{azema, JY2}, the study of random times has always been based on four $\mathbb{F}$ optional processes, which are all defined as projections. Let $\rho$ be a random time. We denote by
$$
\dcb
Z^\rho&\mbox{the optional projection of $\ind_{[0,\rho)}$;}\\
\widetilde{Z}^\rho&\mbox{the optional projection of $\ind_{[0,\rho]}$;}\\
A^\rho&\mbox{the optional dual projection of $\ind_{\{0<\rho\}}\ind_{[\rho,\infty)}$;}\\
a^\rho&\mbox{the predictable dual projection of $\ind_{\{0<\rho\}}\ind_{[\rho,\infty)}$.}
\dce
$$
The process $Z^\rho$ is called the Azema supermartingale of $\rho$. The process $\widetilde{Z}^\rho$ will be called the strong supermartingale of $\rho$. The process $\widetilde{Z}^\rho$ is not cadlag in general. But this process, or more specifically the random set $\mathtt{C}^\rho:=\{\widetilde{Z}^\rho=1\}$, carries the most of knowledge about the random time $\rho$ that we can obtain with an $\mathbb{F}$ optional process. The set $\mathtt{C}^\rho=\{\widetilde{Z}^\rho=1\}$ will be called the optional shadow of $\rho$.  With respect to $\mathtt{C}^\rho$, for any $t\in\mathbb{R}_{+}$, we introduce$$
\dcb
G^\rho_{t}=\sup\{s\leq t: s\in\mathtt{C}^\rho\},&&g^\rho_{t}=\sup\{s< t: s\in\mathtt{C}^\rho\},\\
D^\rho_{t}=\sup\{s> t: s\in\mathtt{C}^\rho\},&& d^\rho_{t}=\sup\{s\geq t: s\in\mathtt{C}^\rho\}.
\dce
$$
Below, we will call $G^\rho_{t}$ the last sojourn times and $D^\rho_{t}$ the entrance times of $\mathtt{C}^\rho$. We set $M^\rho = Z^\rho + A^\rho$ and $m^\rho = Z^\rho + a^\rho$. The processes $M$ and $m$ are $\mathtt{BMO}$-martingales (cf. \cite[Definition 10.6]{HWY}, \cite[Chapitre XX, n$^\circ$74]{DM3} and \cite[Lemme(5.17)]{Jeulin80}).

We have the following properties.
\ebe
\item
According to \cite[Lemme(4.3)]{Jeulin80} \cite[Chapitre XX, n$^\circ$15]{DM3}, $\mathtt{C}^\rho$ is a (random) closed set and it is the biggest optional set contained in $[0,\rho]$. 

\item
The two processes $\widetilde{Z}^\rho$ and $Z^\rho$ are linked by the following relations.$$
\dcb
\widetilde{Z}^\rho_{0}=1, \ Z^\rho_{0}=\mathbb{P}[0<\rho\ |\mathcal{F}_{0}],\\
\widetilde{Z}^\rho = Z^\rho + \Delta A^\rho = Z^\rho_{-} + \Delta M^\rho\ \mbox{ on $(0,\infty)$,}\\
\widetilde{Z}^\rho_{-} = Z^\rho_{-}={^p}\!(\widetilde{Z}^\rho) \ \mbox{ on $(0,\infty)$,}\\
\widetilde{Z}^\rho_{+} = Z^\rho.
\dce
$$

\item
The Azema supermartingale $Z^\rho$ and the strong supermartingale $\widetilde{Z}^\rho$ satisfy the relations below (cf. \cite[Chapitre XX, n$^\circ$30]{DM3} about the relative martingale). For two stopping times $S\leq T$,
$$
\dcb
1 - Z^\rho_{S}
=
\mathbb{E}[(1-Z^\rho_{T})\ind_{\{G_{T}\leq S\}}\ |\mathcal{F}_{S}],\\

1 - \widetilde{Z}^\rho_{S}
=
\mathbb{E}[(1-\widetilde{Z}^\rho_{T})\ind_{\{G_{T}< S\}}\ |\mathcal{F}_{S}],
\dce
$$
on $\{S<\infty\}$.

\item
The set $\{\widetilde{Z}^\rho>Z^\rho\}=\lri 0_{\{1>Z^\rho_{0}\}} \rri\cup\{\Delta A^\rho>0\}$ is an optional thin set, i.e., it is the union of graphs of a sequence of stopping times (cf. \cite[Appendix du chapitre IV, n$^\circ$117]{DM3} and \cite[Definition 3.19]{HWY}). 
\item
$\{Z^\rho_{-}=1\}\subset \mathtt{C}^\rho$ (cf. \cite[Lemme(4.3)]{Jeulin80} or \cite{AMY}) or $\{Z^\rho_{-}=1 \mbox{ or } Z^\rho=1\}\subset \mathtt{C}^\rho$.

\item
$\ind_{\{Z^\rho_{-}=1\}}\Delta M^\rho = \ind_{\{Z^\rho_{-}=1\}}(\widetilde{Z}^\rho-Z^\rho_{-})=0$.

\item\label{ZtildeZ}
Using the identity $\widetilde{Z}^\rho = Z^\rho + \Delta A^\rho$, we see that, for any $b\in(0,1)$, the set $\{Z<b\}\cap\{\widetilde{Z}=1\}$ has no accumulation point.

\item
We introduce, for $n\in\mathbb{N}^*$, $R_{n}:=\inf\{t\geq 0: Z^\rho_{t}\leq \frac{1}{n}\}$. Then, $\rho\in\cup_{n=1}^\infty[0,R_{n}]$ almost surely. Let $R=\sup_{n\geq 1}R_{n}$. Then, $Z^\rho_{s}=0$ for all $s\in[R,\infty)$.
\dbe

\

\section{Honest time}\label{honesttime}

The notion of honest time has been widely studied in the literature See, for example, \cite{azema, barlow, Jeulin80, JY2, MSW}. We will take the following lemma as definition.

\bl
A $\mathcal{F}_{\infty}$ measurable random time $\rho$ is honest, if and only if $\widetilde{Z}^\rho_{\rho}=1$ on the set $\{\rho<\infty\}$, if and only if there exists an optional set $\mathtt{O}$ such that $\rho$ coincides with the $\sup\!\mathtt{O}$ on $\{\rho<\infty\}$. 
\el

\subsection{The first properties of a honest time}\label{basicprop}

For honest time $\rho$, we have the additional properties.
\ebe
\item[$.$]
$\rho=\sup\mathtt{C}^\rho$ on $\{\rho<\infty\}$.

\item[$.$]
Outside of $\mathtt{C}^\rho$, $\widetilde{Z}^\rho = {Z}^\rho$ (cf. \cite[Chapitre XX, n$^\circ$16]{DM3} or \cite[Lemme(5.2)]{Jeulin80}).

\item[$.$]
For a finite random time $\rho$ to be the end of a predictable set, it is necessary and sufficient that $Z^\rho_{\rho-}=1$ on $\{\rho<\infty\}$ (\cite[Proposition 3]{JY2}). In this case the random measure $da^\rho$ is supported on $\{Z^\rho_{-}=1\}$ and $\rho
=
\sup\{Z_{-}=1\}$. See \cite[Chapitre XX]{DM3}.
\dbe

\subsection{Supports of $dA^\rho$ for a honest time}

As in \cite[Chapitre XX]{DM3} we denote by $\mathbf{G}^\rho$ the points in $\mathtt{C}^\rho$ which are isolated at right, i.e., $$
\mathbf{G}^\rho
=
\{t\in\mathtt{C}^\rho: t<D^\rho_{t}\}.
$$
We have the following lemma according to \cite[Theorem 3.3]{AMY} or \cite[Chapitre XX, n$^\circ$13]{DM3}.

\bl\label{exactsupport}
For a honest time $\rho$, the closed support of $dA^\rho$ is the closure of $\mathbf{G}^\rho$. In particular, $\mathtt{C}^\rho$ is a support of the random measure $dA^\rho$.
\el

Consider now the set $\{\Delta A^\rho>0\}$ of the atoms of the random measure $dA^\rho$. We have the following corollary.

\bcor
For a honest time $\rho$, $\{\Delta A^\rho>0\}=(0,\infty)\cap\{\widetilde{Z}^\rho>Z^\rho\}\subset \mathtt{C}^\rho$. 
\ecor

We will denote $\mathbf{G}^\rho_{o}=(0,\infty)\cap\{\widetilde{Z}^\rho>Z^\rho\}=\{\Delta A^\rho>0\}$ (a notation of \cite[Chapitre XX n$^\circ$8]{DM3}).

\

\

\

\

\

\section{The running maximum processes}\label{running}

For any function $u$ defined on $\mathbb{R}_{+}$, we call $$
u^*_{t}=\sup_{0\leq s\leq t}u_{s}, \ \ t\geq 0,
$$
the \textit{running maximum function} of $u$ (or \textit{running maximum process}, if $u$ is a stochastic process). We call any $t\in\{u=u^*\}$ a \textit{maximum time} of $u$. We call $\rho=\sup\{u=u^*\}$ the \textit{last maximum time}.

\subsection{The measure associated with a running maximum function}

It is clear that $u^*$ is a non decreasing function. The following result about the support of the measure $du^*$ is well-known.

\bl\label{supportZ=1}
Suppose that $u$ is cadlag. The measure $du^*$ is carried by $\{u=u^*\}$.
\el

Because of its fundamental role in this paper, let us say why this result is true. For any $s\geq 0$, let $\tau_{s}=\inf\{v> s: u^*_{v}>u^*_{s}\}$. As in \cite{azema} or \cite[Chapitre XX, n$^\circ$13]{DM3}, we introduce the left support set $
\Gamma(u^*)=\{t>0: \forall 0<s<t, u^*_{s}<u^*_{t}\}.
$
Its complementary is given by $$
\Gamma(u^*)^c=\{t>0: \exists 0<s<t, u^*_{s}=u^*_{t}\}
\subset 
\cup_{s\in\mathbb{Q}_{+}}(s,\tau_{s})+\cup_{s\in\mathbb{Q}_{+}:\Delta_{\tau_{s}}u^*=0}[\tau_{s}],
$$
which has a null $du^*$ measure. This shows that $\Gamma(u^*)$ is a support set of $du^*$. Now, we note that, for $t\in\Gamma(u^*)$, either $u_{t}=u^*_{t}$ or there exists a sequence $(s_{k})_{k\geq 1}$ such that $s_{k}\uparrow t$ and $u_{s_{k}}\rightarrow u^*_{t}$. We have therefore
$$
\Gamma(u^*)\subset \{u=u^*\ \mbox{ or }\ u_{-}=u^*\}.
$$
We note that $$
\{u=u^*\ \mbox{ or }\ u_{-}=u^*\}\setminus \{u=u^*\}= \{u< u_{-}=u^*\}.
$$
The set $\{u< u_{-}=u^*\}$ is a denumerable set, and, on this set, $\Delta u^*=0$. This means that $\{u< u_{-}=u^*\}$ has null $du^*$ measure. We conclude that $\{u=u^*\}$ is a support of the measure $d{u^*}$.

As a corollary, we have the next lemma.

\bl\label{U/U}
Let $U$ be a semimartingale with $U_{0}=1$. We have$$
\frac{U}{U^*} = 1+\frac{1}{U^*_{-}} {_{\centerdot}}U - \frac{1}{U^*_{-}} {_{\centerdot}}U^*.
$$
\el

\textbf{Proof.}
The process $\frac{1}{U^*}$ is non increasing and the random measure $d(\frac{1}{U^*})$ is supported on $\{U=U^*\}$ according to Lemma \ref{supportZ=1}. By the integration by parts formula, we write$$
0 = U^* d(\frac{1}{U^*}) + \frac{1}{U^*_{-}} dU^*
\ 
\mbox{ and }\
d(\frac{U}{U^*}) = \frac{1}{U^*_{-}} dU + U d(\frac{1}{U^*}).
$$
These two formulas combined with the support $\{U=U^*\}$ of the random measure $d(\frac{1}{U^*})$ implies $$
\frac{U}{U^*} = 1+ \frac{1}{U^*_{-}} {_{\centerdot}}U - \frac{1}{U^*_{-}} {_{\centerdot}}U^*.\ \ok
$$

\subsection{Skorokhod problem}
One of the reasons of the use of running maximum processes in the study of stochastic processes is its essential role in the resolution of the famous Skorokhod problem.

\bd\textbf{Skorokhod problem}
Let $X$ be a cadlag process. The Skorokhod problem consists to find a process $Y$ such that
\ebe
\item[(i)]
$Y$ is cadlag non decreasing;
\item[(ii)]
$X+Y\geq 0$;
\item[(iii)]
$\int_{\mathbb{R}_{+}}\ind_{\{X_{s}+Y_{s}>0\}}dY_{s}=0$.
\dbe
\ed 

It is well-known that the Skorokhod problem has a solution (cf. \cite{Kons}).

\bl\label{skorokhod}\textbf{Skorokhod Lemma}
The Skorokhod problem has a unique solution given by $Y=(-X)^*\vee 0$.
\el

\subsection{Running maximum process transformation}

When we transform a process, we transform its running maximum process. In some situation, this transformation can be explicitly written down. We begin with a curious result, used in the proof of Theorem \ref{NSforU/U}, on the running maximum processes of stochastic exponentials.

\bl\label{gEg}
Let $w$ and $\gamma$ be two semimartingales null at the origin. Then, $\gamma$ is the running maximum process of $w$ and $\gamma-w\leq 1$, if and only if $1+\mathcal{E}(\gamma)_{-}{_{\centerdot}}w\geq 0$ and $\mathcal{E}(\gamma)$ is the running maximum process of $1+\mathcal{E}(\gamma)_{-}{_{\centerdot}}w$.
\el

\textbf{Proof.}
Suppose that $\gamma$ is the running maximum process of $w$ and $\gamma-w\leq 1$. Note that $\mathcal{E}(\gamma)$ is a non decreasing process whose associated random measure has the same support as the random measure $d\gamma$, contained in $\{w=\gamma\}$. We write$$
\dcb
\mathcal{E}(\gamma)_{-}{_{\centerdot}}w
&=&
\mathcal{E}(\gamma)w - w{_{\centerdot}}\mathcal{E}(\gamma)
=
\mathcal{E}(\gamma)w - \gamma{_{\centerdot}}\mathcal{E}(\gamma)\\
&=&
\mathcal{E}(\gamma)w - \gamma\mathcal{E}(\gamma) + \mathcal{E}(\gamma)_{-}{_{\centerdot}}\gamma
=
\mathcal{E}(\gamma)(w - \gamma) + \mathcal{E}(\gamma) - 1.
\dce
$$
From this identity, we deduce that $1+\mathcal{E}(\gamma)_{-}{_{\centerdot}}w\geq 0$. Denote by $X$ the running maximum process of $1+\mathcal{E}(\gamma)_{-}{_{\centerdot}}w$. For $0\leq s<t$, we have$$
1+\mathcal{E}(\gamma)_{-}{_{\centerdot}}w_{s}\leq \mathcal{E}(\gamma)_{s}\leq \mathcal{E}(\gamma)_{t}.
$$
We prove that $X\leq \mathcal{E}(\gamma)$. On the other hand, let $t'=\inf\{v:\mathcal{E}(\gamma)_{v}=\mathcal{E}(\gamma)_{t}\}$ ($\leq t$). Then, necessarily $t'=\inf\{v:\gamma_{v}=\gamma_{t}\}$ and $t'$ belongs to the left support set $\Gamma(\gamma)$ (cf. the explication of Lemma \ref{supportZ=1} or \cite[Chapitre XX, n$^\circ$13]{DM3} for the definition). Either, $w_{t'}=\gamma_{t'}$ so that $$
\mathcal{E}(\gamma)_{t}=\mathcal{E}(\gamma)_{t'}
=\mathcal{E}(\gamma)_{t'}(w_{t'}-\gamma_{t'})+\mathcal{E}(\gamma)_{t'}
= 
1+\mathcal{E}(\gamma){_{\centerdot}}w_{t'} \leq X_{t},
$$
or $w_{t'}<\gamma_{t'}, w_{t'-}=\gamma_{t'-}$ so that $\mathcal{E}(\gamma)_{t'}=\mathcal{E}(\gamma)_{t'-}$ and $$
\mathcal{E}(\gamma)_{t}=\mathcal{E}(\gamma)_{t'-}
=\mathcal{E}(\gamma)_{t'-}(w_{t'-}-\gamma_{t'-})+\mathcal{E}(\gamma)_{t'-}
= 
1+\mathcal{E}(\gamma)_{-}{_{\centerdot}}w_{t'-} \leq X_{t}.
$$
The about analysis shows that $\mathcal{E}(\gamma)$ is the running maximum of $1+\mathcal{E}(\gamma)_{-}{_{\centerdot}}w$.

Suppose now that $1+\mathcal{E}(\gamma)_{-}{_{\centerdot}}w\geq 0$ and $\mathcal{E}(\gamma)$ is the running maximum of $1+\mathcal{E}(\gamma)_{-}{_{\centerdot}}w$. Let $U$ denote $1+\mathcal{E}(\gamma)_{-}{_{\centerdot}}w$. By Lemma \ref{U/U}, we have$$
\frac{U}{U^*} = 1+ \frac{1}{U^*_{-}} {_{\centerdot}}U - \frac{1}{U^*_{-}} {_{\centerdot}}U^*
=1+w-\gamma.
$$
Rewrite the above equation in the form$$
1- \frac{U}{U^*} = -w +\gamma.
$$
As $U\geq 0$, we conclude that $-w +\gamma\leq 1$. Note that $d\gamma$ has the same left increasing points as $d\mathcal{E}(\gamma)=dU^*$, which are contained in $\{U=U^*\}$. As $1- \frac{U}{U^*}\geq 0$, Lemma \ref{skorokhod} is applicable and we conclude that $\gamma$ is the running maximum process of $w$. \ \ok

\brem
The above lemma can be applied in the following form. \texttt{"}If $U$ is a non negative local martingale with $U_{0}=1$, if $w=\frac{1}{U^*_{-}}{_{\centerdot}}U$ and $\gamma=\frac{1}{U^*}{_{\centerdot}}U^*$, then $\gamma$ is the running maximum of $w$ with $\gamma - w\leq 1$.\texttt{"} \ \ok
\erem

We know that the multiplication of a semimartingale by a positive number does not change the maximum times. This property can be extended to the multiplications by some increasing processes. This extension will be useful in subsection \ref{aboutuniqueness}.

\bpro\label{supmultip}
Let $X$ be a non negative semimartingale. Let $T$ be a stopping time such that $X_{T-}=X^*_{T-}$ on the set $\{T<\infty\}$. Let $v$ be a non decreasing process such that 
\ebe
\item[.]
its associated random measure $dv$ is absolutely continuous with respect to the compensator of the process $\ind_{\{T>0\}}\ind_{[T,\infty)}$, which is therefore supported on the set $\{X_{-}=X^*_{-}\}$, and
\item[.]
$\Delta v\Delta X\ind_{\{\Delta X<0\}}\equiv 0$.
\dbe 
Then, if we define $Y=e^vX$, we have $Y^*=e^vX^*$. 
\epro

\textbf{Proof.} 
$(_{\bullet})$ If $t\in\{X=X^*\}\cap[0,T)$, we have, for $0\leq s\leq t$, $$
Y_{t}e^{-v_{t}}=X_{t}=X^*_{t}\geq X_{s}=Y_{s}e^{-v_{s}},
$$ 
so that $Y_{t}\geq Y_{s}e^{v_{t}-v_{s}}\geq Y_{s}$, i.e., $Y^*_{t}=Y_{t}=X^*_{t}e^{v_{t}}$.

For $t\in\{X<X^*\}\cap[0,T)$, let $$
t'=\sup\{u\leq t: X_{u}=X^*_{t} \ \mbox{ or } X_{u-}=X^*_{t}\}.
$$
Notice that $\{u\leq t: X_{u}=X^*_{t} \ \mbox{ or } X_{u-}=X^*_{t}\}$ is a closed set. We have either $X_{t'}=X^*_{t}$ or $X_{t'}<X_{t'-}=X^*_{t}$, that we denote by writing $X^*_{t}=X_{t'\pm}=Y_{t'\pm}e^{-v_{t'\pm}}$. For $s< t'$, we have $$
Y_{t'\pm}e^{-v_{t'\pm}}=X^*_{t}\geq X_{s}=Y_{s}e^{-v_{s}}.
$$
Hence, $Y_{t'\pm}\geq Y_{s}e^{v_{t'\pm}-v_{s}}\geq Y_{s}$. On the other hand, for $s\in[t',t]$, if $X_{t'}=X^*_{t}$, we have $$
Y_{t'}=e^{v_{t'}}X_{t'}= e^{v_{t'}}X^*_{t}\geq e^{v_{t'}}X_{s} =e^{v_{s}}X_{s} = Y_{s},
$$
because $dv$ is supported on $\{X_{-}=X^*_{-}\}$, which implies that $$
Y^*_{t}=Y_{t'}=e^{v_{t'}}X_{t'}=e^{v_{t}}X^*_{t}.
$$ In the same way,  if $X_{t'}<X_{t'-}=X^*_{t}$, we have $$
Y_{t'-}=e^{v_{t'-}}X_{t'-}= e^{v_{t'-}}X^*_{t}= e^{v_{t'}}X^*_{t} \geq e^{v_{t'}}X_{s}=e^{v_{s}}X_{s} = Y_{s}.
$$ 
Here, $\Delta_{t'} v=0$ because $\Delta_{t'} X<0$. We again conclude that $$
Y^*_{t}=Y_{t'-}=e^{v_{t'-}}X_{t'-}=e^{v_{t}}X^*_{t}. 
$$

$(_{\bullet\bullet})$ Consider $Y^*$ at $T<\infty$. If $X_{T}=X^*_{T}$, we have $$
Y_{T}=e^{v_{T}}X_{T} = e^{v_{T}}X^*_{T}\geq e^{v_{s}}X_{s} = Y_{s},
$$
for all $0\leq s\leq T$, i.e., $Y^*_{T}=Y_{T}=e^{v_{T}}X^*_{T}$. If $X_{T}<X_{T-}=X^*_{T}$, we have $\Delta_{T}v=0$ and $$
Y_{T-}=e^{v_{T-}}X_{T-} = e^{v_{T-}}X^*_{T}= e^{v_{T}}X^*_{T}\geq e^{v_{s}}X_{s} = Y_{s},
$$
for all $0\leq s\leq T$, i.e., $Y^*_{T}=Y_{T-}=e^{v_{T}}X^*_{T}$.

$(_{\bullet\bullet\bullet})$ At the end, consider $Y^*$ at $(T,\infty)$. The process $v$ is stopped at $T$. The relation $Y^*=e^{v}X^*$ on $(T,\infty)$ is straightforward. \ \ok

\

\

\

\section{Martingale maximal representation}\label{DoobMR}

We study now the {martingale maximal representation} (cf. condition \ref{MMR}).

\subsection{The process $(U^*)^{-1}_{-}{_{\centerdot}}U^*$}

The goal of this paper is to find all the random times which satisfy the martingale maximal representation, and, for such a random time, to find a local martingale $U\in\mathcal{M}^+_{1,0}$ which realizes the martingale maximal representation. The first thing to do in this research  is to find a good characteristic of random times which links random times to the martingale maximal representation. In this regard, the direct investigation on $U$ appears unsuccessful. On the contrary, we will see later that the process $(U^*)^{-1}_{-}{_{\centerdot}}U^*$ is a good characteristic. For the moment, we present some of the properties of $(U^*)^{-1}_{-}{_{\centerdot}}U^*$.

Let $\rho$ be a random time. We consider $Z^\rho, A^\rho, a^\rho$, etc., associated with $\rho$. Actually, these notations will be used exclusively for this purpose. We will drop out the superscript $\rho$ from the notations and write $Z,A,a$, etc., for $Z^\rho, A^\rho, a^\rho$, etc.. We denote by $\mathcal{M}^+_{1,0}$ the family of all non negative local martingales starting from 1 vanishing at infinity.  

Notice that an immediate consequence of the martingale maximal representation of $\rho$ is $\mathbb{P}[0<\rho<\infty]=1$. Another immediate consequence is the following proposition.

\bpro\label{bbb}
If $Z$ satisfies the martingale maximal representation with a $U\in\mathcal{M}^{+}_{1,0}$, then $dU^*$ is supported on $\{Z=1\}$ and we have$$
dZ = \frac{1}{U^*_{-}} dU - \frac{1}{U^*_{-}} dU^*.
$$
\epro

\textbf{Proof.}
It is the consequence of Lemma \ref{supportZ=1} and Lemma \ref{U/U}. \ \ok

Because of the above proposition, we have the following relations connecting $\frac{1}{U^*_{-}}{_{\centerdot}} U^*$ with $A=A^\rho, a=a^\rho$. (When $B$ is a non decreasing process, we denote by respectively $B^d$ and $B^c$ its purely jump part and its continuous part.)

\bpro\label{-A+U}
If $Z$ satisfies the martingale maximal representation with a $U\in\mathcal{M}^{+}_{1,0}$, then
$- A + \frac{1}{U^*_{-}}{_{\centerdot}} U^*$ as well as $- a + \frac{1}{U^*_{-}}{_{\centerdot}} U^*$ are local martingales. The predictable dual projection of $\frac{1}{U^*_{-}}{_{\centerdot}} U^*$ is $a$, while the predictable dual projection of $-A^d + \frac{1}{U^*_{-}}{_{\centerdot}} U^{*,d}$ is $A^c - \frac{1}{U^*_{-}}{_{\centerdot}} U^{*,c}$.  
\epro

\textbf{Proof.} The first two assertions are consequences of the fact $Z+A$, $Z+a$  and $Z+\frac{1}{U^*_{-}}{_{\centerdot}} U^*$ are all local martingales. The other assertions are the consequences of \cite[Corollary 5.31]{HWY}.\ \ok

\

Now we look at the jumps of $(U^*)^{-1}_{-}{_{\centerdot}}U^*$. We prove that the process $(U^*)^{-1}_{-}{_{\centerdot}}U^*$ can not jump everywhere.

\bpro\label{DU=0}
If $Z$ satisfies the martingale maximal representation with a $U\in\mathcal{M}^{+}_{1,0}$, we have $[A,U^*]\equiv 0$, and hence $\Delta_{\rho}U^*\ind_{\{\rho<\infty\}}\equiv 0$.
\epro

\textbf{Proof.} Clearly $[A,U^*]$ depends only on the jump processes $\Delta A$ and $\Delta U^*$. The first assertion is true, because $\Delta A>0$ happens on the set $\{Z<1\}$, while $\Delta U^*>0$ happens on the set $\{Z=1\}$ according to Lemma \ref{supportZ=1}. For any stopping time $T$, for any number $b>0$, we have$$
\dcb
\mathbb{E}[\ind_{\{0<T<\infty, 0<\Delta_{T}U^*\leq b\}}\ind_{\{T=\rho\}}\Delta_{T}U^*]
&=&
\mathbb{E}[\ind_{\{0<T<\infty, 0<\Delta_{T}U^*\leq b\}}\Delta_{T}U^*\ \mathbb{E}[\ind_{\{T=\rho\}}\ |\mathcal{F}_{T}]\ ]\\
&=&
\mathbb{E}[\ind_{\{0<T<\infty, 0<\Delta_{T}U^*\leq b\}}\Delta_{T}U^*\ \Delta_{T}A\ ] = 0.
\dce
$$
This proves the second assertion of the proposition. \ \ok

\bpro
Set $\mathbf{D}=\{t\in\mathtt{C}: g_{t}<t\}$. Suppose the martingale maximal representation $Z=\frac{U}{U^*}$ with a $U\in\mathcal{M}^+_{1,0}$. Then, $[Z,U^*]=[M,U^*]$ is a non decreasing purely jumping process and the random measure $d[Z,U^*]$ is supported on the set $\{Z_{-}<1=Z\}\subset\mathbf{D}$. 
\epro

\textbf{Proof.} Let $t>0$ such that $\Delta_{t}U^*>0$. Then, $
U_{t-}\leq U^*_{t-}<U^*_{t}=U_{t}.
$
In order to have $\Delta_{t}Z\neq 0$, the only choice is to have $U_{t-}< U^*_{t-}$. In this case, necessarily $\Delta_{t}Z>0$, i.e., $t\in\{Z_{-}<1=Z\}$. By the property \ref{ZtildeZ} of subsection \ref{basic-elements}, $\{Z_{-}<1=Z\}\subset\mathbf{D}$.\ \ok

\subsection{The case of a honest time $\rho$}

When $\rho$ is in addition a honest time, we can say more about the jumps $(U^*)^{-1}_{-}\Delta U^*$.

\bpro\label{DADGmF=0}
Suppose that $\rho$ is honest and $Z$ satisfies the martingale maximal representation with a $U\in\mathcal{M}^{+}_{1,0}$.  
\ebe
\item
The process $\ind_{\{Z_{-}<1\}}{_{\centerdot}}(-A^d + \frac{1}{U^*_{-}}{_{\centerdot}} U^{*,d})$ is a purely jumping process (a sum of jumps). 
\item
For any $s> 0$, $$
\mathbb{E}[-\Delta_{D_{s}}A + \frac{1}{U^*_{D_{s}-}}\Delta_{D_{s}}U^*\ |\mathcal{F}_{D_{s}-}] = 0,
$$ 
on $\{s<D_{s}<\infty\}$.
\dbe
\epro

\textbf{Proof.} Recall that, when $\rho$ is honest, $dA$ is supported on $\{\widetilde{Z}=1\}$. The martingale maximal representation implies that $dU^*$ is supported on $\{{Z}=1\}$. For the first assertion, it is because $\{\widetilde{Z}=1, Z<1\}$ and $\{Z\neq Z_{-}\}$ are denumerable sets so that $\{\widetilde{Z}=1, Z_{-}<1\}$ and $\{Z=1, Z_{-}<1\}$ are denumerable. For the second assertion, it is because the process $-A + \frac{1}{U^*_{-}}{_{\centerdot}} U^{*}$ on $(s,D_{s}]$ is a local martingale and a step process with only one possible jump at $D_{s}<\infty$ (which is a property related to \cite[Chapitre VI n$^\circ$60]{DM2}).\ \ok

\brem
The jumps $\Delta A$ is also linked with $U$. We have $U^*\Delta A = (U^* - U)\ind_{\mathbf{G}_{o}}$. If $\rho$ is the end of a predictable set, we have $U^*\Delta A = \ind_{\{Z_{-}=1>Z\}}(U^* - U)\ind_{\mathbf{G}_{o}}= \ind_{\{Z_{-}=1>Z\}}(U_{-} - U)= -\ind_{\{U_{-}=U^*>U\}}\Delta U$. \ \ok
\erem

\

\subsection{Exemple}\label{compensated-poisson}

The study of the jumps $(U^*)^{-1}_{-}\Delta U^*$ leads to the following result.

\bpro
The martingale maximal representation is not satisfied by all random time.
\epro

Here is a counter-example. Let $N$ be a unit Poisson process and let $\widetilde{N}_{t}=N_{t}-t, t\geq 0$, be the compensated Poisson process. Then, $W=\mathcal{E}(\widetilde{N})$ is in $\mathcal{M}^+_{1,0}$. The set $\{W=W^*\}$ is almost surely a finite set (so that a closed set). Hence, if $\rho=\sup\{W=W^*\}$, we have $\mathtt{C}^\rho=\mathbf{G}^\rho=\{W=W^*\}$. By Lemma \ref{exactsupport}, for $s> 0$, $\Delta_{D_{s}}A>0$ on the set $\{s<D_{s}<\infty\}$. Consequently, no local martingale $U\in\mathcal{M}^+_{1,0}$ can exist to satisfy the conditions in Proposition \ref{DU=0} and in Proposition \ref{DADGmF=0}:$$
\Delta_{D_{s}}A\Delta_{D_{s}}U^*\ind_{\{s<D_{s}<\infty\}}=0\ \mbox{ and }\ 
\mathbb{E}[-\Delta_{D_{s}}A+\frac{1}{U^*_{D_{s}-}}\Delta_{D_{s}}U^*\ |\mathcal{F}_{D_{s}-}]\ind_{\{s<D_{s}<\infty\}}=0.
$$ 
We conclude that $Z^\rho$ does not satisfy the martingale maximal representation.

\

\subsection{Construction of martingale maximal representation}

We have seen different properties of the process $(U^*)^{-1}_{-}{_{\centerdot}}U^*$ and a counter-example. In particular, by Proposition \ref{-A+U}, the martingale maximal representation can be realized for a random time $\rho$, only if the optional dual projection $A$ of $\ind_{\{\rho>0\}}\ind_{[\rho,\infty)}$ has a \texttt{"}martingale compensator\texttt{"} which lives on $\{Z=1\}$.  In this subsection we will prove that the existence of such a \texttt{"}martingale compensator\texttt{"} for $A$ is also a sufficient condition for the martingale maximal representation.

\bethe\label{NSforU/U}
Suppose that $Z_{0}=1,Z_{\infty}=0$. Then, $Z$ possesses a martingale maximal representation if and only if there exists a process $\gamma$ such that 
\ebe
\item[$.$]
$\gamma$ is non decreasing and $\gamma_{0}=0$; 
\item[$.$] 
the random measure $d\gamma$ is supported on $\{Z=1\}$; 
\item[$.$]
$-A+\gamma$ is a local martingale. 
\dbe
In this case, $Z$ has the martingale maximal representation $Z=\frac{U}{U^*}$, where $U$ is the unique solution of $U = 1+ U^*{_{\centerdot}}(Z+\gamma)$ which is in $\mathcal{M}^+_{1,0}$.

\ethe

\textbf{Proof.} 
\textbf{Part I.}
Suppose the existence of the process $\gamma$. Let $W$ be the local martingale which is the solution of the stochastic differential equation (cf. \cite[Chapter V, Theorem 7]{protter}) $$
W = 1+ W^*_{-}{_{\centerdot}}(Z+\gamma).
$$
Then, writing $$
- \frac{1}{W^*_{-}}{_{\centerdot}}W + \gamma = -(Z-Z_{0} + \gamma) +\gamma = 1-Z,
$$
by Lemma \ref{skorokhod}, we see that $\gamma$ is the running maximum of $\frac{1}{W^*}{_{\centerdot}}W$. By Lemma \ref{gEg}, $
\mathcal{E}(\gamma)
$
is the running maximum of $U:=1+\mathcal{E}(\gamma)_{-}\frac{1}{W^*_{-}}{_{\centerdot}}W$, which is non negative. Now, we compute $$
\dcb
\frac{U}{U^*}
&=&
 1+\frac{1}{U^*_{-}} {_{\centerdot}}U - \frac{1}{U^*_{-}} {_{\centerdot}}U^*
=
 1+\frac{1}{W^*_{-}} {_{\centerdot}}W - \gamma
=
1+ Z-1 = Z.
\dce
$$
Notice that $Z_{\infty}=0$ and $U$ is a non negative local martingale. Necessarily $U\in\mathcal{M}^+_{1,0}$, i.e., $Z$ satisfies the martingale maximal representation. 

Notice that $U_{0}=W_{0}=1$ and $$
\frac{1}{U^*_{-}} {_{\centerdot}}U=\frac{1}{W^*_{-}} {_{\centerdot}}W
= Z-1+\gamma.
$$
By the uniqueness of the stochastic differential equation, we have $U=W$.

\textbf{Part II.}
Suppose now that $Z$ satisfies the martingale maximal representation with a $U\in\mathcal{M}^+_{1,0}$. According to Proposition \ref{bbb} and Proposition \ref{-A+U} the process $\gamma =\frac{1}{U^*_{-}} {_{\centerdot}}U^*$ satisfies the conditions of the theorem.\ \ok

Let us turn back for a moment to the source \cite{AP, K2015, NY} of our discussion and give a different proof of their main result in application of Theorem \ref{NSforU/U}. 

\bpro
Let $\rho$ be the end of a predictable set with $\mathbb{P}[0<\rho<\infty]=1$. Then, $Z$ has the martingale maximal representation with a $U\in\mathcal{M}^+_{1,0}$ such that $U^*$ is continuous, if and only if $\rho$ avoids any predictable stopping time.
\epro

\textbf{Proof.}
We check that, if $\rho$ is the end of a predictable set and avoids the predictable stopping times, the process $\gamma = a$ is continuous and satisfies the conditions of Theorem \ref{NSforU/U} (notably the set $\{Z=1\}$ being a support of $da$). This proves the \texttt{"}if\texttt{"} part of the proposition.

Conversely, if $Z=\frac{U}{U^*}$ for a $U\in\mathcal{M}^+_{1,0}$ such that $U^*$ is continuous, the continuity of $\frac{1}{U^*_{-}} {_{\centerdot}}U^*$ and Proposition \ref{bbb} imply that $a = \frac{1}{U^*_{-}} {_{\centerdot}}U^*$ and consequently $\rho$ avoids any predictable stopping time.  \ \ok

\subsection{About the uniqueness of the martingale maximal representation}\label{aboutuniqueness}

We consider the uniqueness of the martingale maximal representation constructed in Theorem \ref{NSforU/U}. In this section we use the notations introduced in Theorem \ref{NSforU/U}. Suppose that $\gamma$ is a process which satisfies the conditions in Theorem \ref{NSforU/U}. The idea is to show, in modifying $\gamma$, that it may exist different \texttt{"}martingale compensators\texttt{"} of $A$, which live all on $\{Z=1\}$. We can therefore have different martingale maximal representations.

\bpro\label{modif-predictable}
Let $T$ be a predictable stopping time. Then, for any non negative random variable $\xi$ such that $$
\mathbb{E}[-\Delta_{T}A + \xi\ind_{\{\Delta_{T}\gamma>0\}}\ |\mathcal{F}_{T-}]\ind_{\{T<\infty\}}=0,
$$
the process $\widehat{\gamma} := \gamma - \Delta_{T}\gamma\ind_{[T,\infty)} + \xi\ind_{\{\Delta_{T}\gamma>0\}}\ind_{[T,\infty)}$ satisfies the conditions in Theorem \ref{NSforU/U}.
\epro

\textbf{Proof.}
Clearly $\gamma - \Delta_{T}\gamma\ind_{[T,\infty)}$ and $\xi\ind_{\{\Delta_{T}\gamma>0\}}\ind_{[T,\infty)}$ are non decreasing processes. As the support of $d\gamma$ is contained in $\{Z=1\}$, when $\Delta_{T}\gamma>0$, $\lri T\rri\subset\{Z=1\}$. We deduce that $d\widehat{\gamma}$ is supported on $\{Z=1\}$. We write$$
\dcb
-A+\widehat{\gamma}
&=&
-A + \Delta_{T}A\ind_{[T,\infty)} - \Delta_{T}A\ind_{[T,\infty)}
+\gamma - \Delta_{T}\gamma\ind_{[T,\infty)} + \xi\ind_{\{\Delta_{T}\gamma>0\}}\ind_{[T,\infty)}\\

&=&
(-A+\gamma) - \Delta_{T}(-A+\gamma)\ind_{[T,\infty)} + (- \Delta_{T}A+ \xi\ind_{\{\Delta_{T}\gamma>0\}})\ind_{[T,\infty)}.
\dce
$$
We know that $(-A+\gamma)$ is a local martingale. As $T$ is predictable, $\Delta_{T}(-A+\gamma)\ind_{[T,\infty)}$ also is a local martingale. Finally, the condition on $\xi$ implies that $(- \Delta_{T}A+ \xi\ind_{\{\Delta_{T}\gamma>0\}})\ind_{[T,\infty)}$ is again a local martingale.\ \ok

\bpro
Let $T$ be a totally inaccessible stopping time such that $\lri T\rri\subset \{Z_{-}=Z=1\}$. For any non negative random variable $\xi\leq \Delta_{T}\gamma$, let $v,v'$ be respectively the predictable dual projections of the processes$$
\Delta_{T}\gamma\ind_{[T,\infty)}\ \mbox{ and }\ \xi\ind_{\{\Delta_{T}\gamma>0\}}\ind_{[T,\infty)}.
$$
Then, the process $$
\widehat{\gamma} := \gamma - (\Delta_{T}\gamma\ind_{[T,\infty)}-v) + (\xi\ind_{\{\Delta_{T}\gamma>0\}}\ind_{[T,\infty)} - v')
$$ 
satisfies the conditions in Theorem \ref{NSforU/U}.
\epro

\textbf{Proof.}
The process $\gamma - \Delta_{T}\gamma\ind_{[T,\infty)}$ is non decreasing. As $\xi\leq \Delta_{T}\gamma$, the process $v-v'$ is non decreasing. Clearly, the process $\xi\ind_{\{\Delta_{T}\gamma>0\}}\ind_{[T,\infty)}$ is non decreasing. We conclude that $\widehat{\gamma}$ is non decreasing. 

Since $\lri T\rri\subset \{Z_{-}=Z=1\}$, the random measures $dv, dv'$ are supported on $\{Z_{-}=1\}$. As $T$ is totally inaccessible, $v,v'$ are continuous. As $\{Z=1\}\Delta \{Z_{-}=1\}$ is denumerable, we prove that the random measures $dv, dv'$, and hence $d\widehat{\gamma}$, are supported on $\{Z=1\}$.

The process $
-A+\widehat{\gamma}
=
-A+\gamma + (\widehat{\gamma}-\gamma)
$
is a local martingale, because $-A+\gamma$ as well as $(\widehat{\gamma}-\gamma)$ are local martingales. \ \ok

The next proposition follows fundamentally the same idea, but has a proof of different style.

\bpro
Suppose that $Z$ satisfies the martingale maximal representation $Z=\frac{U}{U^*}$ with a $U$ in $\mathcal{M}^+_{1,0}$. Let $T$ be a totally inaccessible stopping time such that $Z_{T-}=1$. If $\Delta_{T}U>0$ on $\{T<\infty\}$, let $v$ be the compensator of the process $\frac{1}{U_{T-}}\Delta_{T}U\ind_{[T,\infty)}$ and define$$
U'= U \frac{e^v}{1+U^{-1}_{T-}\Delta_{T}U\ind_{[T,\infty)}}.
$$
Then, $U'$ is a local martingale in $\mathcal{M}^+_{1,0}$. We have $$
U^*=e^{-v}(1+U^{-1}_{T-}\Delta_{T}U\ind_{[T,\infty)}) U'^*\ \mbox{ and }\ Z=\frac{U'}{U'^*}.
$$ 
\epro

\textbf{Proof.} 
\textbf{Part I.}
We note that $U'$ is well-defined, because $\Delta_{T}U>0$ and $U_{T-}=U^*_{T-}\geq 1$ on $\{T<\infty\}$. Let $b>0$ be a real number. We have $b+U=(b+1)\mathcal{E}(W)$ where $W=(b+U_{-})^{-1}{_{\centerdot}}U$. Let $v''=\frac{U_{-}}{b+U_{-}}{_{\centerdot}}v$ so that $v''$ is the compensator of the process $\frac{1}{b+U_{T-}}\Delta_{T}U\ind_{[T,\infty)}$. Consider the local martingales$$
W''=\frac{1}{b+U_{T-}}\Delta_{T}U\ind_{[T,\infty)} - v'' =\Delta_{T}W\ind_{[T,\infty)} - v'' \ \mbox{ and }\ W'=W -(\Delta_{T}W\ind_{[T,\infty)} - v'')=W-W''.
$$
We have $[W',W'']\equiv 0$ so that (cf. \cite[Problems 9.17]{HWY})$$
b+U=(b+1)\mathcal{E}(W) = (b+1)\mathcal{E}(W')\mathcal{E}(W'')
=
(b+1)\mathcal{E}(W')e^{-v''}(1+\Delta_{T}W\ind_{[T,\infty)}).
$$
Hence $$
(b+1)\mathcal{E}(W')
=
(b+U)\mathcal{E}(W'')^{-1}
=
(b+U)\frac{e^{v''}}{1+\Delta_{T}W\ind_{[T,\infty)}}
=
(b+U)\frac{e^{U_{-}(b+U_{-})^{-1}{_{\centerdot}}v}}{1+(b+U_{T-})^{-1}\Delta_{T}U\ind_{[T,\infty)}}
$$
is a local martingale. Let $S$ be any stopping time such that $U^S$ is a uniformly integrable martingale and $v_{S}$ is bounded. We have$$
b+1=(b+1)\mathbb{E}[\mathcal{E}(W')_{S}]
=
\mathbb{E}[(b+U_{S})\frac{e^{U_{-}(b+U_{-})^{-1}{_{\centerdot}}v_{S}}}{1+(b+U_{T-})^{-1}\Delta_{T}U\ind_{\{T\leq S\}}}].
$$
Let $b\downarrow 0$. By the dominated convergence theorem, 
$$
1=
\mathbb{E}[U_{S}\ \frac{e^{v_{S}}}{1+U_{T-}^{-1}\Delta_{T}U\ind_{\{T\leq S\}}}].
$$
This being true for all the above considered stopping time $S$, we conclude that $U'=U\ \frac{e^{v}}{1+U_{T-}^{-1}\Delta_{T}U\ind_{[T,\infty)}}$ is a local martingale.

\textbf{Part II.} We will repeat partially the computations in Proposition \ref{supmultip}.
Let us study the relation between $U^*$ and $U'^*$. As $Z_{T-}=1$, the random measure $dv$ is supported on $\{Z_{-}=1\}$. As $v$ is continuous, $dv$ is also supported on $\{Z=1\}$. Let $0\leq s\leq t$. If $t\in\{U=U^*\}\cap(0,T)$, we have $$
U'_{t}e^{-v_{t}}=U_{t}=U^*_{t}\geq U_{s}=U'_{s}e^{-v_{s}},
$$ 
so that $U'_{t}\geq U'_{s}e^{v_{t}-v_{s}}\geq U'_{s}$, i.e., $U'^*_{t}=U'_{t}=U^*_{t}e^{v_{t}}$. For $t\in\{U<U^*\}\cap(0,T)$, let $$
t'=\sup\{u\leq t: U_{u}=U^*_{t} \ \mbox{ or } U_{u-}=U^*_{t}\}.
$$
Notice that $\mathtt{F}=\{u\leq t: U_{u}=U^*_{t} \ \mbox{ or } U_{u-}=U^*_{t}\}$ is a closed set. We have either $U_{t'}=U^*_{t}$ or $U_{t'-}=U^*_{t}$ (that we denote by $U^*_{t}=U_{t'\pm}$). Therefore, $$
U^*_{t}=U_{t'\pm}=U'_{t'\pm}e^{-v_{t'\pm}}=U'_{t'\pm}e^{-v_{t}},
$$
because $v$ is continuous and $dv$ is supported on $\{Z=1\}$. On the other hand, for $s\leq t$, $$
U^*_{t}\geq U_{s}=U'_{s}e^{-v_{s}},
$$
so that $U'_{t'\pm}\geq U'_{s}e^{v_{t}-v_{s}}\geq U'_{s}$, i.e., $U'^*_{t}=U'_{t'\pm}=U^*_{t}e^{v_{t}}$. From this relation, we see that $$
Z = \frac{U}{U^*} = \frac{U'e^{-v}}{U'^*e^{-v}} = \frac{U'}{U'^*} \ \mbox{ on }[0,T).
$$
To deal with $Z$ on the interval $[T,\infty)$, let us use the notion $W''$ even when $b=0$, i.e., $
W''=\frac{1}{U_{T-}}\Delta_{T}U\ind_{[T,\infty)} - v.
$ 
We note that $W''$ is a process stopped at $T$ and $e^{-v}(1+U^{-1}_{T-}\Delta_{T}U\ind_{[T,\infty)})=\mathcal{E}(W'')$. As $Z_{T-}=1$ and $\Delta_{T}U>0$ on $\{T<\infty\}$, we have necessarily $Z_{T}=1$. Hence,  $U_{T}=U^*_{T}$ and $U_{T-}=U^*_{T-}$ or equivalently $U'^*_{T-}=U'_{T-}=U'_{T}$ because $\Delta_{T}U'=0$ (as we can check by a direct computation). Hence, for $t>T$,$$
U^*_{t}=\sup_{0\leq s\leq t}U_{s} =\sup_{T\leq s\leq t}U_{s}
=\mathcal{E}(W'')_{T}\sup_{T\leq s\leq t}U'_{s}
=\mathcal{E}(W'')_{T}\sup_{0\leq s\leq t}U'_{s}.
$$
We conclude $Z=\frac{U}{U^*}=\frac{\mathcal{E}(W'')_{T}U'}{\mathcal{E}(W'')_{T}U'^*}=\frac{U'}{U'^*}$ on $[T,\infty)$. As $Z_{\infty}=0$ and $v_{\infty}<\infty$, $U'\in\mathcal{M}^+_{1,0}$. And, combining the above properties, we also see that $U^*=\mathcal{E}(W'')U'^*$. 
\ \ok

\textbf{Example.}
We finish the paper by an example. Let $S$ be a finite totally inaccessible stopping time. Consider the random time $\rho=S$. We have $$
A = \ind_{[S,\infty)}\
\mbox{ and }\
Z=\ind_{[0,S)}.
$$
Let $\gamma$ be the compensator of $\ind_{[S,\infty)}$. Then, $-A+\gamma$ is a local martingale and the random measure $d\gamma$ is carried on $[0,S)=\{Z=1\}$. Hence, by Theorem \ref{NSforU/U}, the martingale maximal representation holds with $U$ determined by the equation $$
U=1+U^*_{-}{_{\centerdot}}(Z+\gamma)=1+U^*_{-}{_{\centerdot}}(-\ind_{[S,\infty)}+\gamma).
$$
Let $U=e^\gamma\ind_{[0,S)}$. Then, we have $U^*=e^\gamma$ and$$
\dcb
1+U^*_{-}{_{\centerdot}}(Z+\gamma)
&=&
1+U^*(Z+\gamma) -1 - (Z+\gamma)_{-}{_{\centerdot}}U^*
=
e^\gamma(Z+\gamma) - (Z+\gamma)_{-}{_{\centerdot}}e^\gamma\\
&=&
e^\gamma Z+e^\gamma \gamma - (e^\gamma-1) - \gamma{_{\centerdot}}e^\gamma\\
&=&
e^\gamma Z+e^\gamma \gamma - (e^\gamma-1) - \gamma e^\gamma + e^\gamma {_{\centerdot}}\gamma\\
&=&
e^\gamma Z - (e^\gamma-1)   + (e^\gamma-1) = e^\gamma Z = U,
\dce
$$
i.e., $U$ is the unique solution of the equation of Theorem \ref{NSforU/U}, which is a local martingale in $\mathcal{M}^+_{1,0}$. Effectively, we check directly $Z=\ind_{[0,S)} = \frac{U}{U^*}.$ 

The question is now if one can construct a different \texttt{"}martingale compensator\texttt{"} of $A$ living on $\{Z=1\}$. Suppose that $S$ is the first jump time of a unit Poisson process so that $\gamma_{t}=t\wedge S$ for $t\geq 0$. Let $S'$ be also the first jump time of a unit Poisson process. Suppose $\mathbb{P}[S= S']=0$. For any bounded predictable process $H$, we have$$
\mathbb{E}[H_{S}]=\mathbb{E}[\int_{0}^SH_{s}ds]\ \mbox{ and }\
\mathbb{E}[H_{S'}\ind_{\{S'< S\}}]=\mathbb{E}[H_{S'}\ind_{\{S'\leq S\}}]=\mathbb{E}[\int_{0}^{S\wedge S'}H_{s}ds],
$$
so that$$
\mathbb{E}[H_{S}] - \mathbb{E}[H_{S'}\ind_{\{S'< S\}}] - \mathbb{E}[\int_{S\wedge S'}^SH_{s}ds]=0.
$$
Or equivalently, $A - \ind_{\{S'< S\}}\ind_{[S',\infty)} - \ind_{(S'\wedge S, S]}.\boldsymbol \lambda$ is a local martingale, where $\boldsymbol \lambda$ denotes the identity map $\boldsymbol \lambda(t)=t, t\geq 0$. The process $\gamma':=\ind_{\{S'< S\}}\ind_{[S',\infty)} + \ind_{(S'\wedge S, S]}.\boldsymbol \lambda$ is non decreasing and $d\gamma'$ is supported on $[0,S)=\{Z=1\}$. We have got a different \texttt{"}martingale compensator\texttt{"} of $A$ living on $\{Z=1\}$.

\end{document}